\documentstyle{article}
\newfont{\petiot}{cmfib8}

\def\C{\mbox{l\hspace{-.47em}C}}

\def\Cn{{\bf C}^{n}} 
\def\TC{T{\bf C}^{n}}
\def\TeC{{T}^{*} {\bf C}^{n}}
\def\sjk{\sum_{j,k=l+1}^{l+m}}
\def\sj{\sum_{j=l+1}^{l+m}}
\def\Tenm{{T}_{N}^{*}M}
\def\Tnm{{T}_{N}M}
\def\doj{\frac{\partial}{\partial {x}_{j}}}
\def\TcM{{T}^{c}M}
\def\Cka{{C}^{(k, \alpha)}, \ k \geq 2, 0 < \alpha < 1}
\def\Cn{{\bf C}^{n}}
\def\Tesm{{T}_{S}^{*}M}
\def\Tsm{{T}_{S}M}
\def\TemC{{T}_{M}^{*} {\bf C}^{n}}

\def\TmC{{T}_{M} {\bf C}^{n}}

\bigskip

\begin{document}
\baselineskip=0.5cm

$\:$
\bigskip
\bigskip
\bigskip
\bigskip

\begin{center}
{\bf  GLOBAL MINIMALITY OF GENERIC MANIFOLDS AND}

{\bf HOLOMORPHIC EXTENDIBILITY OF CR FUNCTIONS}

\bigskip
\bigskip

J. Merker

\bigskip

{\small
Summer address:\\
All\'ee du noyer\\
F-25170 CHAUCENNE\\
FRANCE. \\
\ \\
   Universit\'e Paris VI \\
   ANALYSE COMPLEXE ET GEOMETRIE \\
     U.R.A. 213 du CNRS \\
   Tour 45-46 5\`e \'etage -- Bo\^{\i}te 172 \\
    4, place Jussieu \\
  75252 PARIS CEDEX 05 \\
}

\bigskip
\bigskip
\bigskip
\noindent
{\bf Introduction.}
\end{center}

Let $M$ be a smooth generic submanifold of $\Cn$. 
Several authors have studied the property of CR functions 
on $M$ to extend locally to manifolds with boundary 
attached to $M$ and holomorphically to generic wedges 
with edge M ({\em cf.} \cite{BP}, \cite{TU1}, \cite{TU2}).
In a recent work (\cite{TU3}), Tumanov has showed that CR-extendibility of CR 
functions on M propagates along 
curves that run in complex tangential directions to M.
His main result appears 
as a natural generalization of results by 
Tr\'epreau on propagation of singularities of CR functions (\cite{TR2}).
Indeed, Theorem 5.1 in \cite{TU3} states that the direction of CR-extendibility
 moves parallelly with respect to a certain differential geometric 
partial connection in a quotient bundle of the normal bundle to $M$,
and this variation is dual to the one introduced by Tr\'epreau, 
according to Proposition 7.3 in \cite{TU3}.

In this paper we give a new
and simplified presentation of 
the connection introduced in Tumanov's work.
Let $M$ be a real manifold and $N$ a submanifold of $M$,
$K$ a subbundle of $TM$ with the property that $K {|}_{N} \subset TN$.
Then by means of the Lie bracket, we can define a $K$-partial
connection on the normal bundle of $N$ in $M$ (Proposition 1.1).
In general, the parallel translation associated with 
that  partial connection will be induced by the flow of 
$K$-tangent sections of $TM$ (Proposition 1.2).
When $M$ is a generic submanifold of $\Cn$ containing a CR submanifold $S$
with the same CR dimension we recover in section 2 the ${T}^{c}S$-partial connection 
constructed by Tumanov in \cite{TU3}.

Recall that the {\it CR-orbit} of a point $z \in M$ is  the set of 
points that can  be reached by piecewise smooth integral curves of 
complex tangent vector fields.
We then say that $M$ is {\it globally minimal} at a point 
$z \in M$ if the CR-orbit of $z$ contains a  neighborhood
of $z$ in $M$.
Using previous results, we show that  vector space
generated by 
the directions of CR-extendibility  of CR functions on $M$ exchanges by 
the induced composed flow between two points in a same CR-orbit (Lemma 3.5).
As an application, we prove the main result of this paper, conjectured by Tr\'epreau
in \cite{TR2} :
 {\it for wedge extendibility of CR functions
 to hold at every point in the CR-orbit
of $z \in M$ it is sufficient that $M$ be globally minimal at $z$}
(Theorem 3.4).
Up till now we can only conjecture the converse (for a local result, see \cite{BR90}).

I wish to thank J.-M. Tr\'epreau for helpful
critical and simplifying remarks.

{\bf Remark :} \ After this paper was completed, we have received a preprint by 
B. J\"{o}ricke {\em Deformation of CR-manifolds, minimal points and CR-manifolds
with the microlocal analytic extension property},
which contains also a proof of Theorem 3.4 and  Theorem 3.6.
Our proof seems quite different since we obtain these results
relying on Tumanov's propagation theorems,
the generic manifold $M$ being fixed, whereas
 B. J\"{o}ricke
works with conic perturbations of the base manifold so as to produce 
minimal points.

\bigskip 
\begin{center}
{\bf  \S 1. Partial connections associated with  a system of
vector fields.}
\end{center}

\medskip

Let $M$ be a real differentiable manifold of  class  ${C}^{2}$
of dimension $n$
and $H \rightarrow M$ a $r$-dimensional vector bundle over $ M$.
Recall that a connection $\nabla$
on the bundle $H \rightarrow M$ is a bilinear mapping 
which assigns to each pair of a vector field $X$ with domain $U$  
and a section  $\eta$ of $H$ over $U$ a section ${\nabla}_{X}  \eta$
of $H$ over $U$ and satisfy
 $${\nabla}_{\phi X} =   \phi {\nabla}_{X},  \ \ \ 
 {\nabla}_{X} (\phi \eta)  =  \phi {\nabla}_{X} \eta  + (X \phi) \eta ,
 \ \ \ \ \ \ \ \ \ \ \phi \in {C}^{1}(M,{\bf R}). $$

When the covariant derivative ${\nabla}_{X} \eta$ can only be defined for 
vectors $X$ that belong to a subbundle $K$ of $TM$, we call the connection 
$\nabla$ a {\it K-partial} connection (cf. \cite{TU3}).

If $N$ is a submanifold of $M$, let ${T}_{N} M$ be the {\it 
normal bundle of N in M}, {\em i.e.}
$${T}_{N} M \ \ = \ \ {TM|}_{N}/ TN .$$

\noindent
{\bf \sc  Proposition 1.1.} \ {\it Let M be a real manifold of class ${C}^{2}$,
$N \subset M$ a submanifold of class ${C}^{2}$ too and let K be a ${C}^{1}$ 
subbundle of $TM$ with the property that 
${K|}_{N} \subset TN$.
Then there exists a natural K-partial connection $\nabla$
 on the bundle ${T}_{N} M$ 
which is defined as follows.
If $x \in N$,$X \in K [ xÊ]$ and $\eta$ is a local section of ${T}_{N} M$
over a neighborhood of $x$,  then take 
$$
{\nabla}_{X} \eta \ \ = \ \ [\tilde{X}, \tilde{Y}](x) \ \ mod \ {T}_{x} N
$$
where $\tilde{X}$ is  a ${C}^{1}$  local section of $K$ extending $X$ 
 and $\tilde{Y}$ is a 
lifting of $\eta$ in $TM$ in a neighborhood of x.}

\medskip 
{\sc Proof.}
We first check that the definition is independent of the lifting 
$\tilde{Y}$. 
In fact, when $\tilde{Y}$ is  tangent to $N$, as $\tilde{X}$ is tangent to $N$ too,
the Lie bracket $[\tilde{X}, \tilde{Y}]$ remains tangent to $N$ hence is zero
in the quotient bundle. 

Next we have to check that the definition of $\nabla$ is independent of the chosen section 
$\tilde{X}$ or, to rephrase, that if $\tilde{X} (x) = 0$ then 
$[\tilde{X}, \tilde{Y}](x)$ belongs to 
${T}_{x}N$.
Since $K$ is a fiber bundle we can write 
$$\tilde{X} \ \ =  \ \ \sum_{j=1}^{r} {f}_{j} {\tilde{X}}_{j}  
\ \  \ \ \ \ \ {f}_{j} (0) = 0  \ \ \ \ \ j=1, ..., r $$
where $r = rank  \ K$, ${({\tilde{X}}_{j})}_{j=1,...,r}$ \
 \ is a frame for $K$ near $x$ and
the ${f}_{j}$ are ${C}^{1}$ real valued functions defined near $x$.
Noting that 
$$[f \tilde{X}, \tilde{Y}] \ \equiv \ f[\tilde{X}, \tilde{Y}] -
(\tilde{Y} f) \tilde{X}  
\ \equiv  \ f[\tilde{X}, \tilde{Y}] \ \  \ \ \ \ mod  \ TN
$$
the result follows and the mapping $\nabla$ is 
well-defined. 
Moreover the preceding implies that if $\phi
\in {C}^{1} (M, {\bf R})$
$$
{\nabla}_{\phi X} \eta  \  \equiv
 \ \phi {\nabla}_{X} \eta .
$$
Last, we check that ${\nabla}_{X}$ is a derivation.
Indeed
$$
{\nabla}_{X} (\phi \eta) \ \equiv  \ [\tilde{X}, \phi \tilde{Y}](x) \ \equiv \
(\tilde{X}.\phi) \tilde{Y} + \phi [\tilde{X}, \tilde{Y}] \ \equiv \
(X\phi)  \eta + \phi {\nabla}_{X} \eta
$$
and the proof is complete.

\smallskip
With the connection $\nabla$ it is associated
 the {\it parallel translation} of fibers of
${T}_{N} M$ along smooth curves on the base $N$ that 
run in directions tangent to $K$.
Let $I \ni t$ be a subinterval of ${\bf R}$ and
 $ \gamma : I \rightarrow N $ be a smooth curve with the property that 
$\dot{ \gamma}(t) \in K [ \gamma(t) ]$,
where $\dot{ \gamma} = \frac{d}{dt}  \gamma(t)$.
A curve $ \eta (t) \in {T}_{N}M[ \gamma(t) ]$ is a {\it horizontal lift} of $ \gamma $
if ${\nabla}_{\dot{ \gamma }} \eta = 0$.
Existence and uniqueness of horizontal lifts provide 
 linear 
isomorphisms
$$
{\Phi}_{{t}_{0},t} \ \ \ : {T}_{N}M[ \gamma ({t}_{0})] \ \rightarrow \ 
{T}_{N}M[ \gamma(t)  ]
$$
obtained by moving elements of $K$ along horizontal lifts of 
$\gamma$.

Recall (cf. \cite{SPI}) that the Lie bracket $[\tilde{X}, \tilde{Y}]$ 
is  defined as the 
Lie derivative ${L}_{\tilde{X}} \tilde{Y}$ of $\tilde{Y}$ with respect to 
$\tilde{X}$ 
$$
[ \tilde{X}, \tilde{Y}](x) \ = \ {L}_{\tilde{X}} \tilde{Y} \ = \
 \lim_{h \rightarrow 0} [ \tilde{Y}(x) - d {\tilde{X}}_{-h} 
(\tilde{Y} ({\tilde{X}}_{h}(x))) ]
$$
where ${\tilde{X}}_{t}$
is the local flow on $M$ generated in a neighborhood of $x$ by  the  
vector field
$\tilde{X}$, and $d {\tilde{X}}_{t}$ denotes its differential.
In the assumptions of Proposition 1.1, $\tilde{X}$
is of class ${C}^{1}$ so the mapping 
$x \rightarrow {\tilde{X}}_{t} (x)$ is of class
${C}^{1}$ and the differential is a well-defined continuous 
mapping.
When $\tilde{X}$ is $ K$-tangent its flow (and more generally any piecewise smooth
composition of such flows) stabilizes the tangent bundle
$TN$ of the manifold $N$, hence its differential induces isomorphisms of fibers of 
${T}_{N}M$, which we denote by $d {X}_{t}$. 
Assume moreover that the curve
$\gamma$ is  an integral curve of a ${C}^{1}$ $K$-tangent vector field
$\tilde{X}$, (which cannot be true for  most general smooth curves
$\gamma$ but is sufficient enough for the applications) : 
 $\gamma (0) = x$ 
and $\gamma (t) \ = \ {\tilde{X}}_{t} (x)$.
Then we claim that the mapping
$$
d {X}_{t} \  \ \ : \ \ \ {T}_{N}M [x] \ \rightarrow \ {T}_{N}M [ \gamma (t)]
$$
 provides the parallel translation ${\Phi}_{0,t}$.
Indeed let ${\eta}_{0} \in {T}_{N}M [x]$ and take $\eta (t) = d {X}_{t} ({\eta}_{0})$.
Then by the definition of the partial connection $\nabla$
and the definition of the Lie bracket we have
$$
{\nabla}_{\dot{\gamma}} \eta (t) \  \ \equiv  \ \ 0 .
$$
By uniqueness of solutions of 
linear differential equations of order one it must be that
$$
\eta (t)  \ \ = {\Phi}_{0,t} ( {\eta}_{0}) .
$$

\bigskip
{\sc Proposition 1.2.} \ {\it
Under the hypotheses of Proposition 1.1,
let $\gamma (t) = {X}_{t} ({x}_{1})$
be a smooth (piecewise smooth) integral curve of a K-tangent vector field X (a
finite number of K-tangent vector fields)
running from $x_{1} \in N$ to $x_{2} \in N$.
Then the parallel translation along $\gamma$ associated with the K-partial
connection $\nabla$ 
is induced by the differential of the flow of X (composed flow).
}

\smallskip

In order to give an expression of 
the covariant derivatives induced by the partial
 connection $\nabla$, we
 choose coordinates on $M$, $x=(x', x'') \in {\bf R}^{l} \times {\bf R}^{m}$ such
that the base point corresponds to $x=0$
 and the submanifold $N$ is defined  by the equation $x''=0$.
Let $(x, \eta)= (x',x'',\eta ', \eta '')$ be the canonical coordinates on $TM$, and
$(x', \eta '') \in {\bf R}^{l} \times {\bf R}^{m}$ the associated 
coordinates on $\Tnm$.

If 
$X= \sum_{j=1}^{l+m} {a}_{j} (x)  \ \frac{\partial}{\partial {x}_{j}}$
is a ${C}^{1}$ section of $K$, it must be tangent to N, so
${a}_{j} (x', 0) = 0, \ j=l+1,...,l+m$.
We choose a local section $\eta$ of $\Tnm$ over a neighborhood of $0$ in $N$, in 
fact a section $\tilde{Y}$ of $TM$  of the form
$$\tilde{Y}= \sum_{j=l+1}^{l+m} {\eta}_{j}(x') \  \doj ,$$
Recalling Proposition 1.1
 we have the following expression for the covariant derivative
of $\eta$ in the direction of $X$
$${\nabla}_{X} \eta =
   \sum_{j=l+1}^{l+m} \sum_{k=1}^{l}
    {a}_{k} (x',0) \ \frac{\partial {\eta}_{j}}{\partial {x}_{k}} (x')  \  \doj \ - \
    \sjk {\eta}_{k} (x') \ \frac{\partial {a}_{j}}{\partial {x}_{k}} (x',0) \  \doj .$$
Given an integral curve
$\gamma (t) = (\gamma ' (t), 0)$
 of the field $X$, the equations for the horizontal 
lifts look like
$$X. {\eta}_{j} =
      {\dot{\eta}}_{j}(t)=
          \sum_{k=l+1}^{k=l+m} \frac{\partial {a}_{j}}{\partial {x}_{k}} (\gamma'(t),0) 
          {\eta}_{k} (t) \ \ \ \ \ j=l+1,...,l+m$$
so the curve $(\gamma ' (t), \eta''(t))$
 is the integral curve of the following vector 
field $\check{X}$ on $\Tnm$ 
$$\check{X}(x', \eta'') = \sum_{j=1}^{l} {a}_{j}(x',0) \doj \ + \ 
    \sjk \frac{\partial {a}_{j}}{ \partial {x}_{k}} (x',0) {\eta}_{k}(x') 
\frac{\partial}{ \partial {\eta}_{j}} .$$
Alternately, the partial connection $\nabla$
can be defined by the family of horizontal 
subspaces 
$H(\eta) \subset {T}_{\eta} (\Tnm)$
generated by vectors of the form $\check{X}$.

\medskip
Let us consider the {\it dual connection} 
${\nabla}^{*}$ to the connection $\nabla$ on the dual
bundle $\Tenm$.
Recall that the conormal bundle of 
$N$ in $M$,  $\Tenm$, 
consists of forms in ${T}^{*}M$ that vanish on $TN$.
It has fiber over a point $ x \in N$
$$\Tenm [x] = \{ \phi \in {T}_{x}^{*}M ; \ \phi {|}_{{T}_{x}N} = 0 \} .$$
The dual connection
${\nabla}^{*}$ is defined by the  following relation :
if X is a $K$-tangent vector to $N$ at $x$, $\eta$ is any section of $\Tnm$ 
near  $x$ and $\phi$ is any section of $\Tenm$
$$X< \phi, \eta> \ = \ <{\nabla}_{X}^{*} \phi, \eta > + <\phi, {\nabla}_{X} \eta> .$$
It is easily checked that such a relation defines a 
$K$-partial connection on $\Tenm$.

Along with the coordinates on $\Tnm$ we introduced before we can introduce
the canonical coordinates $(x', \xi '')$ on the conormal bundle $\Tenm$.
These are dual to the coordinates $(x', \eta '')$
for the canonical duality $<,>$ between
$\Tnm$ and $\Tenm$ :
$$< \sj {\xi}_{j} d {x}_{j} \ , \ \sj {\eta}_{j} \doj > \ = 
\ \sj {\xi}_{j} {\eta}_{j} .$$
Using the previous definition of the dual connection 
we can then compute the covariant derivative of a section 
$\sum {\xi}_{j} d {x}_{j} = \phi$ of $\Tenm$.
One easily shows
$${\nabla}^{*}_{X} \phi = \sj ( X. {\xi}_{j} + \sum_{k=l+1}^{l+m} 
            {\xi}_{k} \frac{\partial {a}_{k}}{ \partial {x}_{j}}) \ d x_{j} .$$
Hence, 
under the assumption of Proposition 1.2,
  the parallel translation associated with the connection 
${\nabla}^{*}$ is given by means of the integral curves of the following vector field
on $\Tenm$
$$\hat{X}= \sum_{j=1}^{l} {a}_{j} (x', 0) \doj -
    \sjk \frac{\partial {a}_{j}}{\partial {x}_{k}} (x', 0) {\xi}_{j} 
          \frac{\partial}{\partial {\xi}_{k}} .$$

\smallskip
There is another way of thinking the connection ${\nabla}^{*}$ dual
to the partial connection $\nabla$ which has been considered by 
Tr\'epreau  in \cite{TR2}.

To a general vector field $X$ on $M$ it is 	associated its symbol $\sigma (X)$
which is an invariantly defined function on the cotangent bundle
${T}^{*}M$ of $M$.
To a function $f$ of class ${C}^{1}$ on ${T}^{*}M$ it is associated its hamiltonian
field ${H}_{f}$.

Let ${X}_{j}, \ j=1,...,r$ \  be a local basis of $K$-tangent sections of $TM$.
Let ${\Sigma}_{K}$ be the orthogonal complement of $K$ in 
${T}^{*}M$.
If $X = \sum_{j=1}^{r} {\phi}_{j} {X}_{j} $
is a ${C}^{1}$ section of $K$ we have
$${H}_{\sigma (X)} {|}_{{\Sigma}_{K}} \ = \
   \sum_{j=1}^{r} {\phi}_{j} {H}_{\sigma ({X}_{j})} {|}_{{\Sigma}_{K}} +
   \sum_{j=1}^{r} \sigma ({X}_{j}) {H}_{{\phi}_{j}} {|}_{{\Sigma}_{K}} .$$
Since $\sigma ({X}_{j}) ,  \ j=1,...,r$ is zero on
${\Sigma}_{K}$, we deduce that the  restricted hamiltonian field
$${H}_{\sigma (X)} {|}_{{\Sigma}_{K}} \ = \
\sum_{j=1}^{r} {\phi}_{j} {H}_{\sigma ({X}_{j})} {|}_{{\Sigma}_{K}} $$
depends only on the value of X at the base point and not on the 
chosen section.
If $X$ is tangent to $N$,
${H}_{\sigma (X)}$ when restricted to $\Tenm$ is tangent to $\Tenm$.
Hence we have constructed another vector field
on $\Tenm$ which is in fact the same as the one associated with the 
connection dual to the partial  connection $\nabla$.

Indeed, let as before $(x', \xi'')$ be the canonical
coordinates on the conormal bundle $\Tenm$.
Recall that the hamiltonian field of a function $f=
f(x, \xi)$  just looks like
$$H_{f}= \sum_{j,k=1}^{j,k=l+1}  \frac{\partial f}{\partial {\xi}_{k}} \doj - 
\frac{\partial f}{\partial {x}_{j}} \frac{\partial}{\partial {\xi}_{k}} . $$ 
The symbol of the section
$$X= \sum_{j=1}^{l+m} {a}_{j} \doj \ \ \ \ \ {a}_{j} (x', 0)=0, \ j=l+1,...,l+m$$
of K being $\sigma (X) = \sum {a}_{j} {\xi}_{j}$
we can compute 
$${H}_{\sigma (X)} {|}_{\Tenm} =
   \sum_{j=1}^{l} {a}_{j} (x',0) \doj
    - \sjk \frac{\partial {a}_{j}}{\partial {x}_{k}} (x',0) {\xi}_{j} 
                \frac{\partial}{\partial {\xi}_{k}}$$
and the last expression proves that  ${H}_{\sigma (X)}$
is the same vector field on $\Tenm$ as $\hat{X}$
computed previously,
so the set of  restricted hamiltonian fields 
${H}_{\sigma (X)} {|}_{\Tenm}$
defines the same family of horizontal subspaces
for the partial connection ${\nabla}^{*}$.

The next section is devoted to the application of the preceding results to 
the geometry of CR submanifolds of $\Cn$.

\bigskip
\begin{center}
{\bf \S 2. Application to generic submanifolds of $\Cn$}
\end{center}

\medskip
In this section we apply results 
 of section 1 in the context of  
differential geometry in the complex euclidean space $\Cn$.
Afterwards we check that our definitions recover
those of Tr\'epreau \cite{TR2} and Tumanov \cite{TR3}.

\smallskip
Let $\TC$ be the real tangent bundle of $\Cn$ and $J$
be the standard complex structure operator on $\TC$.
Let $\TeC$ be the bundle of {\it holomorphic } (${\bf C}$-linear) 1-forms
on $\Cn$.
In the canonical coordinates $z = (z_{1},...,z_{n})$ its fiber
over a point $z$ consists of (1,0)-forms
$\omega = \sum_{j=1}^{n} {\zeta}_{j} d {z}_{j}, \ {\zeta}_{j} \in {\bf C}, \
j=1,...,n$.
Then $\TeC$
is a {\it complex} manifold. It can be (and it is usually) identified with the real
dual bundle of $\TC$ introducing the real duality defined by 
$$(\omega, X) \ \in \TeC \times \TC \ \ \ \ \ \ \ \ 
(\omega, X) \ \longmapsto \ Im \ <\omega, X>.$$
In other words we identify real and holomorphic forms
by $Im \ \omega \leftrightarrow \ \omega$.

Now, let $M$ be a real submanifold of $\Cn$.
In this identification, the conormal
bundle $\TemC$
is a subbundle of $\TeC$ and it has fiber spaces
$$\TemC [ z ] \ = \ \{ \omega \in \TeC; \  Im \  \omega  \ {|}_{{T}_{z}M} = 0 \ \} .$$
Hence the bundles $\TmC = \TC {|}_{M} / TM$ and $\TemC$ are in duality by
$$(\omega, X \ mod \ TM) \ \ \longmapsto \ \ {\rm Im} \ <\omega, X> .$$

Assume moreover that $M$ is generic (that is $TM + JTM = \TC {|}_{M}$)
and let ${\Sigma}_{M}$ be the orthogonal complement of the complex
tangent bundle $\TcM$ in the cotangent bundle ${T}^{*}M$.  In the
terminology of linear partial differential equations it is the {\it
characteristic set} (and since $\TcM$ is a fiber bundle, the
characteristic {\it manifold}) of the system of CR vector fields.  It
is easily checked that ${\Sigma}_{M}$ and $TM / \TcM$ are in duality
in the same way.

Since $M$ is CR, ${\Sigma}_{M}$
is a fiber bundle and there is a canonical bundle epimorphism
$$\theta \ : \ \ \TemC \ \rightarrow \ {\Sigma}_{M} ,$$
defined by 
$\theta (\omega) = {\imath}^{*}_{M} \omega$
where ${\imath}_{M} : \ M \ \rightarrow \ {\bf C}^{n}$
is the natural injection.
Since $M$ is generic, $\theta$ is an isomorphism.

On the other hand the complex structure $J$ induces an isomorphism,
still denoted by $J$
$$J \ : \ \ TM / \TcM \ \rightarrow \ \TmC .$$

\medskip
{\sc Lemma 2.1.} \ {\it $\theta$ is the transposed of J, i.e. }
$$(\omega, JX) \ = \ (\theta (\omega), X)$$
{\it for every } $\omega, \ X$.

\medskip
(Indeed $<\omega, JX>=i<\omega,X>$).

\smallskip
From now on we let $S \subset M$ be a CR submanifold of $M$ with the 
property that $CRdim \  S = CRdim \ M$.
Equivalently it is required that ${T}^{c} S = \TcM {|}_{S}$.
By restriction analogous pairs of bundles remain isomorphic when 
$TM / \TcM$
is replaced by $\Tsm$, 
${\Sigma}_{M}$ is replaced by $\Tesm$,
$\TmC$ is replaced by 
${T {\bf C}^{n} |}_{S} \ / \ 
({TM |}_{S} + JTS) = E$,
and $\TemC$ is replaced by
$\TemC \cap i {T}_{N}^{*} {\bf C}^{n} = {E}^{*}$,
but now ${T}^{c}S$-partial connections can be defined by means of the 
isomorphisms $J$ and $\theta$ 
on the two new bundles $E$ and ${E}^{*}$.
Note that the duplication essentialy deals with 
complex differential geometry.

First, the results of the previous section apply
with $K = \TcM = TM \cap JTM$
and $N = S$ and produce a ${T}^{c}S$-partial connection 
$\nabla$ on $\Tsm$
together with the dual connection 
${\nabla}^{*}$ on $\Tesm$.
On the other hand,
the push forward by $J$ of $\nabla$ defines a
${T}^{c}S$-partial connection $\Theta$ on $E$ ;
its action on a section $\vartheta$ of $E$ in the
 direction of a complex tangent
vector $X$ is simply
$${\Theta}_{X} \vartheta \ \ := \ \ J \ {\nabla}_{X} \ ( {J}^{-1} \ \vartheta).$$

Similarly, the pull-back of the
${T}^{c}S$-partial connection 
${\nabla}^{*}$ by $\theta$ defines a ${T}^{c}S$-partial connection 
${\Theta}^{*}$ on ${E}^{*}$, and ${\Theta}^{*}$ is the connection dual to $\Theta$
since $\theta $
is the transposed of $J$ (lemma 2.1).

\smallskip
Recall from section 1 that if $X$ is a section of $\TcM$
then $\hat{X} = {H}_{\sigma (X)} {|}_{\Tesm}$
is tangent to $\Tesm$.
In \cite{TR2}, Tr\'epreau showed that ${E}^{*}$
is a CR manifold, using a lemma which states that given such a vector field 
$\hat{X}$ tangent to $\Tesm$
with horizontal part $X$ complex tangent to $M$, there exists a unique
vector field $\tilde{X}$ complex tangent to ${E}^{*}$ with the same 
horizontal part $X$.
Moreover Tr\'epreau states that
$$\hat{X} \ = \ d \theta ( \tilde{X})$$
Hence we deduce that the ${T}^{c}S$-partial connection
${\Theta}^{*} \ = \ {\theta}^{*}  {\nabla}^{*}$ can alternately
be given, as is originally done in \cite{TR2}, by means of the vector fields
of the form $\tilde{X}$, i.e. horizontal
subspaces of  
${\Theta}^{*}$ are spanned by tangent vectors
to integral curves of $\tilde{X}$.
We then have checked that the parallel translation
in ${E}^{*}$ introduced by Tr\'epreau
with the assumption of Proposition 1.2
 is the 
same as the one associated with the ${T}^{c}S$-partial
 connection ${\Theta}^{*}$ previously defined
starting, as in section 1, with the partial connection
associated with the bundle of complex tangents to $M$, $K= \TcM$.
Moreover, since $\tilde{X}$ is complex tangent to ${E}^{*}$,
we see that
${T}^{c} {E}^{*}$
is the set of horizontal subspaces for the ${T}^{c}S$-partial
 connection ${\Theta}^{*}$ .
This has been noticed in \cite{TU3} Êand will be usefull in the next section
when proving Theorem 3.4.

\bigskip
\begin{center}
{\bf \S 3. Orbits and the extension of CR functions.}
\end{center}

\medskip

\smallskip
In this section, it is assumed that
 $M$ is a generic submanifold of $\Cn$ of smoothness class 
${C}^{2}$,
and we let ${\bf X}$ be the set of ${C}^{1}$ sections over open subsets of $M$
of $\TcM$.
If $z \in M$, the  subset of $M$ consisting of points of $M$ which can be reached by piecewise
${C}^{1}$-smooth integral curves of elements of ${\bf X}$, starting at $z$,
is called the {\it CR-orbit} of $z$, and is 
denoted by $\cal{O}$  [$z$].

If $U$ is an open subset of M, ${\bf X} {|}_{U}$ denotes the set
of elements of ${\bf X}$ restricted to $U$.
It is well-known ({\em cf.} \cite{SU}, \cite{BR90}, \cite{TR2}) that

\vspace{-0.7cm}

\ \ \ \ \ \ \ \ \ \ \ \ \ \ \ \ \ \ \ \ \  \ \ \ \ \ \ \ \ \
\begin{picture}(220,60) 
\put(75,20){$\longrightarrow$}
\put(75,25){lim}
\put(80,10){U}
\put(100,25){$\cal{O}( {\bf X} {|}_{U},z$)}
\end{picture}

\vspace{-0.3cm}
\noindent
where $U$ runs over the open neighborhoods of $z$ in $M$ defines the germ at $z$ of
the unique CR-submanifold  of $M$ 
with the same CR dimension as $M$
of minimal dimension passing through $z$, 
which is called the
{\it local CR-orbit} of $z$ and is
denoted by ${\cal{O}}^{loc}  \ [z]$.
When   considering  ${\cal{O}}^{loc}  \ [z]$ in the following
we shall mean such a submanifold of a neighborhood
of $z$ in $M$, i.e. an actual representative of the germ.
It plays the crucial role in the study of automatic
extendibility of CR functions ({\em cf.} Theorem 3.1 below).

\medskip
 Recall that a smooth complex-valued function on $M$ is called a {\it CR function}
if it is annihilated by every antiholomorphic tangent vector field on $M$.
A continuous function can be thought CR in the sense of distribution theory.
We denote by $CR(M)$ the set of all continous CR functions on $M$.

For completeness we recall definitions from \cite{TR2} and \cite{TU3}.
We say that a manifold $\tilde{M}$ with boundary is {\it attached to M at (m,u)},
$m \in M$, $u \neq 0, u \in \TmC$[$m$] if $b \tilde{M} \cap U = M \cap U$ for 
some neighborhood $U$ of $m$, and $u$ is represented by a vector ${u}_{1} \in {T}_{m} \tilde{M}$
directed inside $\tilde{M}$.

Let $f$ be a CR function on $M$; we say that $f$ is {\it CR-extendible} at 
$(m,u)$ if it extends continuously to be CR 
on some $\tilde{M}$ attached to $M$ at 
$(m,u)$.
When there is a CR submanifold $S$  of $M$ through $m$
 and a manifold $\tilde{M}$ attached
to $M$ at $(m,u)$, $ u \in \TmC$[$m$],
we also say that $\tilde{M}$  is attached to $M$ 
at $(m, \eta)$, if $u$ represents
 $\eta \in {E}_{m}, \ \eta \neq 0$
($E$ is the bundle defined in section 2).
Similarly it makes sense to consider
 CR-extendibility at $(m, \eta)$, $m\in S, \eta
\in {E}_{m}$.
But it should be noted that given 
$\eta \neq 0$ in ${E}_{m}$
does not determine $\tilde{M}$ unambiguously
 unless $S$ is complex

From now on we will require that $M$ belong 
to the class 
$\Cka$.
This regularity assumption can be justified
since it behaves well when proving  the strongest
local results on CR-extendibility
(In fact, it behaves well through the so-called
{\it Bishop equation}, \cite{TU2}, Theorem 1.), and constructing
wedges with ribs and  an egde having such a regularity (cf. \cite{AY}).
Moreover, we need manifolds of class at least 
${C}^{2}$ in order to apply Proposition 1.1.
Since it will be of use in the proof of Theorem 3.4
we recall the following theorem due to Tumanov (\cite{TU2})

\bigskip {\sc Theorem 3.1.} \ ({\small \sc A. E. Tumanov})
 \ {\it
Let M be a generic submanifold of ${\bf C}^{n}, n=p+q$, 
with dim M = $2p+q$, CRdim M = p, and of smoothness class
${C}^{k,\alpha} \ (k \geq 2), 0  < \alpha < 1$.
For every point $z \in M$ there exist r = r(z) = dim
${\cal O}^{loc}_{[z]} - 2  \ CRdimM$ 
manifolds with boundary ${\tilde{M}}_{1},...,{\tilde{M}}_{r}$
attached to M at z, of class ${C}^{(k, \beta)}$
whenever $0 < \beta < \alpha$ such that}

\ \ \ \ \ \ \ \ \ \ (a) {\it Every CR function on M is CR-extendible to
} ${\tilde{M}}_{1},...,{\tilde{M}}_{r}$

\ \ \ \ \ \ \ \ \ \ (b) $\sum_{j=1}^{r}  {T}_{z'} {\tilde{M}}_{j} =
{T}_{z'}M + J{T}_{z'} {\cal O}^{loc}_{[z]} 
$  \ \ \ \ \ \ \ \ \ \ \ \ {\it z' close to z in} ${\cal O}^{loc}_{[z]}$.

\noindent
{\it Moreover the manifold germ ${\cal O}^{loc}_{[z]}$ is of class
${C}^{(k, \beta)}$ whenever $0 < \beta < \alpha$.}

\smallskip
\smallskip
Note that ${\cal O}^{loc}_{[z]}$ is at least of 
class ${C}^{2}$ so it can play the role of $N$ in 
Propositions 1 and 2.  
Using the connections constructed in section 2
 we can reinterpret the main result on propagation of analyticity for
CR functions recently proved by Tumanov.

According to Tumanov (\cite{TU3}, Proposition 7.3),
the connection dual to the one that is constructed during the paper
has the property that its horizontal subspaces
are exactly fibers of the complex tangent bundle ${T}^{c} {E}^{*}$,
hence, concludes Tumanov,
the induced parallel translation need be the same 
as the one  introduced on ${E}^{*}$ by Tr\'epreau.
We have shown in section 2 that our  connection
$\Theta$ has as a dual connection a connection
${\Theta}^{*}$ with the same property;
so  $\Theta = {J}_{*} \nabla$
coincides with the connection constructed by Tumanov.

Proposition 1.2 together with Theorem 5.1 in \cite{TU3} leads to

\bigskip
{\sc Theorem 3.2.} {\it Let M $\subset {\bf C}^{n}$ be a generic manifold and S
$\subset$ M a CR submanifold of M with the property that 
CRdim S = CRdim M.
Let $\gamma$ be a piecewise smooth integral curve of $\TcM$
running from $z' \in S$ to $z'' \in S$  
 and let ${\Phi}_{\gamma}$
be the associated composed flow.
Then for every $\epsilon > 0$, every $\eta ' \in {E}_{z'}$
and every manifold $\tilde{M} '$ attached to M at ($z',\eta '$), there
exists another manifold $\tilde{M} '$ attached to M at 
($z'', \eta ''$), $\eta '' \in {E}_{z''}$ such that}

 \ \ \ \ \ \ \ \ (a) \ \ $| \eta '' - Jd{\Phi}_{\gamma}(z).{J}^{-1} \eta '| \ < \  \epsilon$

 \ \ \ \ \ \ \ \ (b) \ \ {\it if a CR function on M extends to be CR on}
$\tilde{M}'$ {\it it extends to be CR on }
$\tilde{M}''$

 \ \ \ \ \ \ \ \ (c) \ \ {\it if M}, $\tilde{M}'$ {\it belong to} 
${C}^{k,\alpha} \ (k \geq 2), 0 < \gamma < \alpha < 1$
{\it then there exists such a }
${\tilde{M}}'' \in {C}^{ (k, \gamma) }$.

\medskip
Theorem 3.2  shows that the so-called propagation of analyticity for CR functions
 is intrinsically related to the geometry of the base manifold M.
Moreover, it fundamentally means that the study of extendibility for CR 
functions is closely related to the 
study of 
sections of the complex tangent space to $M$.

\smallskip
Following Sussmann (\cite{SU}), we begin with some adapted terminology and recalls.
Let $X \in {\bf X}$ be a local section of ${T}^{c}M$.
The ${C}^{1}$ integral curves $t \rightarrow \gamma (t)$ of $X$ generate 
local diffeomorphisms of $M$ where they are defined (the so-called {\it flow} of $X$)
 which we will denote by 
$z \rightarrow {X}_{t} z$.
Composites of several maps of the form ${X}_{t}$ can produce local diffeomorphisms of 
neighborhoods of points that are {\it far }from each other in a same CR-orbit. 
If $X = ({X}_{1},...,{X}_{m})$ is an element of ${\bf X}^{m}$ such that
for $t = ({t}_{1},...,{t}_{m}) \ \in {\bf R}^{m}$, the map
 $z \rightarrow   {X}_{m,{t}_{m}} \cdots {X}_{1,{t}_{1}} z$
is well defined in a neighborhood of $z$, we will still  denote it for convenience 
by ${X}_{t}$ or $\Phi$ ({\em cf.} Proposition 1.2).

Let ${\Delta}_{\bf X}$ be the {\it distribution spanned by} ${\bf X}$,
i.e. the mapping which to $z \in M$  assigns the linear hull of vectors 
$X(z)$ where $X$ belongs to ${\bf X}$ : it is just the distribution associated 
with the complex tangent bundle of $M$.
We let ${P}_{\bf X}$ denote the smallest distribution which contains ${\Delta}_{\bf X}$
and is invariant under complex-flow 
diffeomorphisms, or for short the smallest {\it ${\bf X}$-invariant} distribution 
which contains  ${\Delta}_{\bf X}$.
Precisely, ${P}_{\bf X}(z)$ is the linear hull of vectors of the form
$d{X}_{t}(v)$ where $v \in {\Delta}_{\bf X}(z')$ and $z = {X}_{t}z'$.
A ${C}^{1}$ distribution $P$ on $M$ has the 
{\it maximal integral manifold property} if 
for every $z \in M$ there exists a submanifold $S$
 of $M$ such that $z \in S$ and for every
$z' \in  S$,  ${T}_{z'} S = P(z')$.
Moreover, S is said to be a {\it maximal integral manifold} of P
if $S$ is an integral manifold of $ P$ such that every connected 
integral manifold of $P$
which intersects $S$ is an open submanifold of $S$.

Then the results of Sussmann, which extend to the ${C}^{2}$ case tell us 
that $\cal{O}$ [$z$] is a (connected) maximal integral submanifold
of ${P}_{\bf X}$ (perharps with a finer topology) and admits a unique differentiable
 structure
making the injection 
$i : \cal{O}$ [$z$] $\rightarrow M$ an immersion of class ${C}^{1}$.

We now introduce the following definitions.

\bigskip
{\sc Definition 3.3.} \ Let $M$ be a generic submanifold of 
$\Cn$ and $z \in M$.
$M$ is called {\it minimal at z} if ${\cal{O}}^{loc} \ [z]$ contains
a  neighborhood of $z$ in $M$.
It is called {\it globally minimal at z} if $\cal{O}$   [$z$]
contains a neighborhood of $z$.

\smallskip
In view of the global results of Sussmann definition 3.3 means that  
the generic manifold  
 $M$ is globally minimal at a point $z$ if and only if
there exist a finite number of
points $z_{l}^{'}, \ l=1,...,d$ in the CR-orbit of $z$ and
 composed flow diffeomorphisms 
${\Phi}^{'}_{l}, \ l=1,...,d$ of a neighborhood of $z_{l}^{'}$ in $M$
on a neighborhood of $z$ in $M$ respectively such that
$${T}_{z}M \ = \sum_{l=1}^{d} d {\Phi}_{l}({z}^{'}_{l}). \  ({T}^{c}_{z_{l}^{'}} M).$$  

We are now able to prove the theorem conjectured by Tr\'epreau in \cite{TR2}
which is the natural generalization of a celebrated theorem of Tumanov (\cite{TU1}).
Here is the substance of this paper.

\bigskip
{\sc Theorem 3.4.} \ \ {\it Let M be a generic submanifold of $\Cn$ of smoothness class
$\Cka$ which is {\rm globally minimal} at a point $z \in M$.
Then for every z' in the CR-orbit of z
 there exists a wedge $\cal{W}$ of edge M at z'
 such that}

\ \ \ ($ \ast$) \ every CR function on M extends 
holomorphically into $\cal{W}$.

\medskip
{\sc Proof.}
We shall make use of the following abuse of langage :
we will say that a CR-function $u$ is CR-extendible in the direction 
$v \in TM / \TcM [z] $
if it is in fact CR-extendible in the direction of $Jv$.
Let us consider the set
$${H}_{z} = Vect \ \ \ \{ v \in {T}_{z}M / {T}^{c}_{z}M ; \ \ u \ is  \ {\rm CR-extendible}
\  at \ (z,v) \ \}$$
and its preimage under the natural surjection $\pi : TM \rightarrow 
TM / \TcM$
$${\hat{H}}_{z} = {\pi}^{-1} \ ( {H}_{z}) \ \subset \ {T}_{z}M .$$

\bigskip
{\sc Lemma 3.5.} \ {\it Let X be a ${C}^{1}$ section of $\TcM$
over a neighborhood of $z \in M$
and let ${\Phi}_{t}$ be the flow of $X$
and $\Phi = {\Phi}_{t}$ for some $t$.
Then, if $v \in {T}_{z}M$, }
$$v \ \in \ {\hat{H}}_{z} \ \  \  \ \ \ \ \ \Longleftrightarrow \ \ \ \ \ \ \ \
d \Phi (z) .v  \ \in \ {\hat{H}}_{\Phi (z)}
 .$$

\medskip
{\sc Proof.} \ Since the statement is a 
symetric and a transitive one we can assume that 
$z$ and $z'$ are so close that $z':= \Phi (z)$ is contained in a CR submanifold $S$
of $M$ with $CRdim S = CRdim M$ which is minimal at $z$ (for instance
take for $S$ the local CR-orbit of $z$)
and such that $z'$ belongs to the boundary of the manifolds
whose existence comes from Theorem 3.1. Hence 
$$ (*) \ \ \ \ \ \ \ \ \ {\hat{H}}_{z'} \ \supset \ {T}_{z'}S .$$
So if $v$ belongs to ${T}_{z}S$ there is nothing to add.
On the other hand, if $\xi= {pr}_{{T}_{S}M}v \neq 0$ we apply the propagation
result Theorem 3.2 
 and obtain that for every $\epsilon > 0$ $u$ is CR-extendible
at $(z', \xi '')$, where $\xi ''$
is $\epsilon$-close in euclidean norm to
 $\xi ' = d \Phi (x). \xi $ ;
 so letting $\epsilon$ decrease to zero,
since every finite-dimensional vector space is closed, 
we have $\xi ' \in {pr}_{{T}_{S}M} \ ({H}_{z'})$.
Because of ($\ast$) the indetermination on the specific representative
of $\xi'$ is removed whence 
$$d \Phi (z).v  \ \in \ {\hat{H}}_{z'}$$
and the lemma is proved. 

\medskip
{\sc end of proof of theorem 3.4.} \ The global lemma 3.5 and the condition of global 
minimality implie immediately that
$${\hat{H}}_{z'} \ \ = \ \ {T}_{z'}M$$
for every $z'$ in the (global) CR-orbit of $z$.
The conclusion follows by the edge-of-the-wedge theorem  and 
the proof is complete. 

Theorem 4.1 admits an obvious generalization which
involves the concept of ${\cal{W}}_{r}$-{\it wedges}.
Recall that a ${\cal{W}}_{r}$-wedge at $z$
with edge $M$ is locally the general intersection of a wedge of edge $M$
at $z$ and a generic manifold containing $M$ as a 
submanifold of codimension $r$.

\bigskip
{\sc Theorem 3.6.} \ \ {\it Let M be a generic submanifold of $\Cn$ of smoothness class
$\Cka$, and let $r= dim \cal{O}$ [$z$] $-2 CRdim \ M$.
Then for every z' in the CR-orbit of z, every $\gamma$ with 
$0< \gamma < \alpha$, there exists a ${\cal{W}}_{r}$-wedge $\cal{W}$ of edge M at z'
and of smoothness class ${C}^{k, \gamma}$ such that}

\ \ \ ($ \ast$) \ every CR function on M extends to be CR on 
 $\cal{W}$.

\noindent
{\it Moreover, the tangent space to $\cal{W}$ at $z'$ spans
${T}_{z'}M + J {T}_{z'} \ \cal{O}$ [$z$].}

\medskip
{\sc Proof}. The same argument runs
in proving that ${\hat{H}}_{z'}$ contains ${T}_{z'}M + J {T}_{z'} \ \cal{O}$ [$z$]
 and the conclusion then follows by the edge-of-the wedge
theorem of Ayrapetyan (\cite{AY}), in the classes $\Cka$ .

\end{document}